\documentclass[10pt]{amsart}
\usepackage{amsthm}
\usepackage[cyr]{aeguill}
\usepackage{dsfont}
\usepackage{fancyhdr}
\usepackage{upgreek}
\usepackage{mathrsfs}
\usepackage{eufrak}
\usepackage{amsmath,amsfonts,amssymb,color,multicol,subeqnarray}\usepackage{graphicx,longtable}
\usepackage[T1]{fontenc}
\usepackage[latin1]{inputenc}  
\usepackage[a4paper]{geometry}
\usepackage[english]{babel}
\theoremstyle{plain}
\newtheorem{teo}{Theorem}
\newtheorem{pro}{Proposition}
\newtheorem{lem}{Lemma}
\newtheorem{cor}{Corollary}
\theoremstyle{definition}
\newtheorem{defi}{Definition}

\theoremstyle{remark}
\newtheorem*{nota}{\bf Remark}

\def\sca #1_#2_#3 {\langle#1,#2\rangle_#3}

\def\m{\mu}
\def\<{\langle}
\def\>{\rangle}

\def\n{\nu}
\def\pii{\mathcal{\pi}}

\def\ep{\epsilon}

\def\pp{\varphi}

\newcommand{\bq}{\begin{equation}}

\newcommand{\eq}{\end{equation}}
\newcommand{\ba}{\begin{eqnarray*}}
\newcommand{\ea}{\end{eqnarray*}}

\newcommand{\dif}[2]{\frac {\partial #1}{\partial #2}}
\newcommand{\diff}[2]{\frac {\partial^2 #1}{\partial {#2}^2}}

\newcommand{\ban}{\begin{eqnarray}}
\newcommand{\ean}{\end{eqnarray}}

\title{Asymptotic Harmonic Analysis on the Space of Square Complex Matrices}
\date{\today  }
\begin{document}
\maketitle
\centerline{M. RABAOUI\footnote{Institut de Math\'ematiques de Jussieu, Universit\'e Pierre et Marie Curie, 175 rue de chevaleret, 75013 Paris. E-mail : rabaoui@math.jussieu.fr}}

\begin{abstract} In this paper, we determine the spherical functions of positive type on the space $V_\infty= M(\infty, \Bbb C)$ relatively to the action of the product group $K_\infty=U(\infty)\times U(\infty)$. The space $V_\infty$ is the inductive limit of the spaces of square complex matrices $V_n=M(n, \Bbb C)$. The group $K_\infty$ is the inductive limit of the product groups $K_n=U(n)\times U(n)$, where $U(n)$ is the unitary group.
\end{abstract}

\bigskip
\bigskip
\bigskip
\section{Introduction} The work that we present here take place within the framework of the infinite dimensional harmonic analysis on the spherical pairs. We consider in this article the spherical pair $(G_\infty, \ K_\infty)$, which is the inductive limit of the sequence of Gelfand pairs $(G_n, \ K_n)$: $$G_n=K_n\ltimes V_n, \ \ K_n=U(n)\times U(n), \ \ V_n=M(n, \Bbb C), $$ $$G_\infty=K_\infty\ltimes V_\infty, \ \ K_\infty=U(\infty)\times U(\infty).$$ Here $V_\infty=M(\infty,\Bbb C)$ is the space of infinite complex matrices having only a finite number of non-zero entries, and $U(\infty)$ is the group of the infinite unitary matrices $(u_{ij})$ with complex coefficients such that $u_{ij}=\delta_{ij}$ for $i+j$ large enough.

Let $\mathfrak{P}$ be the set of $K_\infty$-biinvariant continuous functions of positive type on $G_\infty$ satisfying $\pp(0)=1$. We are interested in the determination of the extreme points $\pp$ of this convex set which are indeed the spherical functions of postif type relatively to the pair $(G_\infty, \ K_\infty)$. The group $G_\infty$ is equipped with the inductive limit topology. The subgroup $K_\infty$ is closed. The homogeneous space $G_\infty/K_\infty$ is the vector space $$V_\infty=\bigcup_{n=1}^\infty V_n.$$ The law of the inductive limit group $G_\infty=K_\infty\ltimes V_\infty$ is given by : $$(u, \ x)(v, \ y)=\big((u_1v_1, \ u_2v_2), \ x+u_1yu_2^*\big), $$ where $$u=(u_1, u_2), \ v=(v_1, v_2)\in K_\infty \ {\rm and} \ x, \ y\in V_\infty.$$ A function $\pp$ on $G_\infty$ which is right invariant under $K_\infty$ does not depend on the variable $u\in K_\infty$. Therefore it is possible to see it like a function on $V_\infty$ : $$\pp(g)=\pp\big((u, x)\big)=\pp_0(x).$$ Moreover, if the function $\pp$ is $K_\infty$-biinvariant then the function $\pp_0$ which is defined on $V_\infty$ is $K_\infty$-invariant, or $U(\infty)$-biinvariant. Also, the function $\pp$ is of positive type on $G_\infty$ if and only if the function $\pp_0$ is of positive type. Let us note by $D_\infty$ the subspace of diagonal matrices in $V_\infty$. An element of $D_\infty$ can be decomposed as diag$(a_1, a_2, \dots)$, with $a_1,a_2,\dots\in\Bbb R$ and $a_j=0$ for $j$ large enough. Any matrix $x\in V_\infty$ can be diagonalised as $$x=u\ {\rm diag}(a_1, a_2, \dots) \ v^* \ \ \big(u, v\in U(\infty)\big).$$ Consequently, any $K_\infty$-invariant function on $V_\infty$ is uniquely determined by its restriction to the subset $D_\infty$.

For a number of spherical pairs, the spherical functions of positive type (the extreme points of $\mathfrak{P}$) have been determined. One can quote, for example, the work of Schoenberg \cite{Sho1} on $\Bbb R^{(\infty)}$ and those of G. Olshanski and A. Vershik \cite{Ol2} on the space $Herm(\infty,\Bbb C)$ of infinite dimensional Hermitian matrices. In these cases and others the spherical functions of positive type are obtained as limits of the ones of classical Gelfand pairs. In our situation, we obtain the spherical functions in this way. The principal result of this article is the following theorem :\\ \begin{paragraph}{\bfseries{Theorem}} The spherical functions of positive type $\pp$ on $V_\infty$, which are $U(\infty)$-biinvariant and satisfying $\pp(0)=1$, are given by : $$\pp\big({\rm diag}(\xi_1, \dots, \xi_n, 0, \dots)\big)=\Pi(\omega,\xi_1)\dots\Pi(\omega,\xi_n), $$ where $$\Pi(\omega,\lambda):=e^{-\frac{1}{4}\gamma\lambda^2}\prod_{k=1}^\infty\frac{1}{1+\frac{1}{4}\alpha_k\lambda^2}, $$ with $$\omega=(\alpha,\gamma),\ \gamma\in\Bbb R_+, \ \alpha_k\in\Bbb R_+ \ {\rm and} \ \sum_{k=1}^\infty \alpha_k <\infty.$$
\end{paragraph}

\bigskip

\section{Spherical functions of positive type on $(G_n,K_n)$} In this section, we will determine the spherical functions of positive type relatively to the pair $(G_n, K_n)$. An explicit formula of these functions was established for the first time in \cite{Ber}. It was also done in \cite{Mea} by a method using the Abel transform. It was also obtained in \cite{Ben}, by means of a contraction, starting from the spherical functions on $SU(n,n)/S\left(U(n)\times U(n)\right)$. The method that we use here is similar to the one used by J. Faraut in \cite{Far4} for the case of finite dimensional Hermitian matrices. This gives a simple and new proof of the result.\\

 We consider the space $V_n$ of complex matrices on which the group $K_n$ acts as follows : $$T(k):x\mapsto k.x=uxv^* \ \ (u, v\in U(n)).$$ Every matrix $x\in V_n$ admits a polar decomposition $$x=u\,{\rm diag(}\lambda_1, \dots, \lambda_n{\rm)}\,v^* \,, \ \ u, v\in U(n), \, \lambda_j\in\Bbb R. $$ Hence any function $f$ which is $K_n$-invariant on $V_n$ depends only on $\lambda={\rm diag(}\lambda_1, \dots, \lambda_n{\rm)}$ : $$f(x)=F(\lambda_1,\dots,\lambda_n),$$ where $F$ is a function defined on $\Bbb R^n$, invariant under $\mathfrak{S}_n\ltimes\{-1,1\}^n$.
 
 \bigskip

 Let us consider on $V_n$ the Euclidean structure defined by $ \<x,y\> = \mathfrak{Re}\,{\rm tr}(xy^*).$ The laplacian associated to this structure is : \ba\Delta=\sum_{i,j=1}^n\Big(\diff{}{(\mathfrak{Re}\,x_{jj})}+\diff{}{(\mathfrak{Im}\,x_{jj})}\Big)
         &+&\sum_{j<k}^n\Big(\diff{}{(\mathfrak{Re}\,x_{jk})}+\diff{}{(\mathfrak{Re}\,x_{kj})}\Big)\\
         &+&\sum_{j<k}^n\Big(\diff{}{(\mathfrak{Im}\,x_{jk})}+\diff{}{(\mathfrak{Im}\,x_{kj})}\Big).\ea
         
\bigskip

The laplacian is invariant under the action of $K_n$ in the following sense: if $f$ is a function of class $\mathscr{C}^2$, then $$\Delta f\big(T(k)\big)=\Delta\Big(f\big(T(k)\big)\Big)\ \ (k\in K_n).$$

Let $f$ be a function of class $\mathscr{C}^2$ on $V_n$ which is invariant under the action of $K_n$: $$f(uxv^*)=f(x)\ \ (k=(u,v)\in K_n).$$ The function $\Delta f$ is also invariant under $K_n$. This leads to the introduction of the operator $L$ defined by : $$\Delta f(x)=LF(\lambda_1,\dots,\lambda_n).$$ The operator $L$ is called the {\it radial part} of the laplacian.

\bigskip

\begin{pro}\label{radprop} {\rm (i)} Let $f$ be a $K_n$-invariant function of class $\mathscr{C}^2$. Then $$\Delta f(x)=LF(\lambda_1,\dots,\lambda_n),$$ where 
\ba LF=\sum_{i=1}^n\Big(\diff{F}{{\lambda_i}}+\frac{1}{\lambda_i}\dif{F}{{\lambda_i}}\Big)
  &+&2\sum_{i<j}\frac{1}{\lambda_i-\lambda_j}\Big(\dif{F}{{\lambda_i}}-\dif{F}{{\lambda_j}}\Big)\\
  &+&2\sum_{i<j}\frac{1}{\lambda_i+\lambda_j}\Big(\dif{F}{{\lambda_i}}+\dif{F}{{\lambda_j}}\Big).\ea {\rm (ii)} The preceding formula can also be written as
$$LF=\frac{1}{D(\lambda)}\,\sum_{i=1}^n\Big(\diff{}{\lambda_i}+\frac{1}{\lambda_i}\dif{}{\lambda_i}\Big)\big(D(\lambda)F(\lambda)\big),$$ where $D$ is given by $$D(\lambda)=\prod_{i<j}(\lambda_i^2-\lambda_j^2).$$
\end{pro} 

\bigskip

In order to prove the preceding proposition, we will use the following result (see \cite{Far4}, Lemma IX-2.2) :

{{\it {\underline{\bf Lemma :}} Let $f$ be a function of class $\mathscr{C}^2$ on an open set $\mathscr{U}$ of a finite dimensional real vector space $\mathscr{V}$. Let $\mathscr{A}$ be an endomorphism of $\mathscr{V}$, $a\in\mathscr{V}$. Let $\ep>0$ be such that, for $|t|<\ep$, $\exp{t\mathscr{A}}.\,a\in\mathscr{U}$. We assume that, for $|t|<\ep$, $$f(\exp{t\mathscr{A}.\,a)}=f(a).$$ 
Then \ba&(Df)_a\left(\mathscr{A}.\,a\right)=0,&\\ &(D^2f)_a\left(\mathscr{A}.\,a,\,\mathscr{A}.\,a\right)+\left(Df\right)_a\left(\mathscr{A}^2.\,a\right)=0.&\ea}}

\bigskip

\begin{paragraph}{\bfseries{Proof of Proposition \ref{radprop}.}}{\rm Let $\mathscr{U}$ be an open set in $V_n$ and $\mathscr{A}$ the endomorphism of $V_n$ defined by : $\mathscr{A}.\,a=Xa+aY^*$ where $X,Y\in V_n$. If the matrices $X,\,Y$ are skewHermitian, then for every $t\in\Bbb R$, the matrices $\exp{tX}$, $\exp{tY}$ are unitary and, for every $a\in\mathscr{U}$,\\ $$f(\exp{tX}a\exp{tY^*})=f(a).\\$$

We deduce from (\cite{Far4}, lemma IX-2.2) that\\ $$\left(Df\right)_a\left(Xa+aY^*\right)=0,$$ $$\left(D^2f\right)_a\left(Xa+aY^*,\,Xa+aY^*\right)+\left(Df\right)_a\left(X^2a+2XaY^*+a\left(Y^*\right)^2\right)=0 .$$
\bigskip

(a) Let us put $X=Y=E_{jk}-E_{kj}$ ($j\ne k$), $a={\rm diag(}a_1,\dots,a_n{\rm )}$. We obtain\\ $$Xa+aY^*=\left(a_k-a_j\right)\left(E_{jk}+E_{kj}\right),$$ $$X^2a+2XaY^*+a\left(Y^*\right)^2=2\left(a_k-a_j\right)\left(E_{jj}-E_{kk}\right),$$ and hence $$\left(a_k-a_j\right)^2\left(D^2f\right)_a\left(E_{jk}+E_{kj},\,E_{jk}+E_{kj}\right)+2\left(a_k-a_j\right)\left(Df\right)_a\left(E_{jj}-E_{kk}\right)=0 ,$$
where
 $$\diff{f}{(\mathfrak{Re}\,x_{jk})}(a)+\diff{f}{(\mathfrak{Re}\,x_{kj})}(a)={2\over{\left(a_j-a_k\right)}}\left(\dif{f}{(\mathfrak{Re}\,x_{jj})}(a)- \dif{f}{(\mathfrak{Re}\,x_{kk})}(a)\right).$$
 \bigskip
 
 (b) Let us put $X=i\left(E_{jk}+E_{kj}\right)$ et $Y=-X$. We get\\
 $$Xa+aY^*=\left(a_j+a_k\right)\left(iE_{jk}+iE_{kj}\right),$$ $$X^2a+2XaY^*+a\left(Y^*\right)^2=-2\left(a_j+a_k\right)\left(E_{jj}+E_{kk}\right),$$ and hence $$\left(a_j+a_k\right)^2\left(D^2f\right)_a\left(iE_{jk}+iE_{kj},\,iE_{jk}+iE_{kj}\right)-2\left(a_j+a_k\right)\left(Df\right)_a\left(E_{jj}+E_{kk}\right)=0 ,$$
where
 $$\diff{f}{(\mathfrak{Im}\,x_{jk})}(a)+\diff{f}{(\mathfrak{Im}\,x_{kj})}(a)={2\over{\left(a_j+a_k\right)}}\left(\dif{f}{(\mathfrak{Re}\,x_{jj})}(a)+ \dif{f}{(\mathfrak{Re}\,x_{kk})}(a)\right).$$
 \bigskip
 
 (c) Let us put $X=iE_{jj}$ and $Y=-X$. We obtain\\
 $$Xa+aY^*=i\,2a_jE_{jj},$$ $$X^2a+2XaY^*+a\left(Y^*\right)^2=-4a_jE_{jj},$$ and hence
$$4a_j^2\left(D^2f\right)_a\left(iE_{jj},\,iE_{jj}\right)-4a_j\left(Df\right)_a\left(E_{jj}\right)=0 ,$$
where
 $$\diff{f}{(\mathfrak{Im}\,x_{jj})}(a)={1\over{a_j}}\dif{f}{(\mathfrak{Re}\,x_{jj})}(a).$$
 
 Finally, \ba\diff{f}{(\mathfrak{Re}\,x_{jj})}(a)&=&\diff{F}{\lambda_j},\\ \diff{f}{(\mathfrak{Im}\,x_{jj})}(a)&=&{1\over{\lambda_j}}\dif{F}{\lambda_j},\\
\diff{f}{(\mathfrak{Re}\,x_{jk})}(a)+\diff{f}{(\mathfrak{Re}\,x_{kj})}(a)&=&{2\over{\left(\lambda_j-\lambda_k\right)}}\left(\dif{F}{\lambda_j}(a)- \dif{F}{\lambda_k}(a)\right),\\
\diff{f}{(\mathfrak{Im}\,x_{jk})}(a)+\diff{f}{(\mathfrak{Im}\,x_{kj})}(a)&=&{2\over{\left(\lambda_j+\lambda_k\right)}}\left(\dif{F}{\lambda_j}(a)+ \dif{F}{\lambda_k}(a)\right).\ea}

This proves (i). In order to prove (ii), we have to use the formula : $$\Delta_0(DF)=\Delta_0F+2(\nabla_0D|\nabla_0F)+\Delta(D),$$ where $\Delta_0$ is the laplacian and $\nabla_0$ the gradiant on $\Bbb R^n$. The polynomial $D$ is harmonic and $$\sum_{j=1}^n{1\over{\lambda_j}}\dif{D}{\lambda_j}=0.$$ We can then conclude that $${1\over D}\Delta_0(DF)+{1\over D}\sum_{j=1}^n{1\over{\lambda_j}}\dif{(DF)}{\lambda_j}=\Delta_0F+ 2{1\over D}(\nabla_0D|\nabla_0F)+\sum_{j=1}^n{1\over{\lambda_j}}\dif{F}{\lambda_j}.$$

Since $${1\over D}\nabla_0D=\nabla_0\log|D|=\sum_{j<k}\frac{1}{\lambda_j-\lambda_k}(e_j-e_k)+\sum_{j<k}\frac{1}{\lambda_j+\lambda_k}(e_j+e_k),$$ where $(e_1,\dots,e_n)$ is the canonical basis of $\Bbb R^n$, we obtain $${1\over D}(\nabla_0D|\nabla_0F)=\sum_{j<k}\frac{1}{\lambda_j-\lambda_k}\left(\dif{F}{\lambda_j}-\dif{F}{\lambda_k}\right)+\sum_{j<k}\frac{1}{\lambda_j+\lambda_k}\left(\dif{F}{\lambda_j}+\dif{F}{\lambda_k}\right).\hspace{1cm}\Box$$
\end{paragraph}

\bigskip

In the preceding Euclidean polar decomposition, the measure $m$ can be written as : \\ $$\alpha(du)\alpha(dv)\prod_{j<k}\big(\lambda_j^2-\lambda_k^2\big)^2\prod_{j=1}^n\lambda_j\,d\lambda_j ,$$ where $\alpha$ is the normalized Haar measure of the unitary group $U(n)$. Moreover, one has the following integration formula :\\

\begin{pro}\label{intprop}{{\rm (see \cite{Far3}, Proposition X.3.4)}} For every integrable function $f$ on $V_n$\\ $$\int_{V_n}f(x)m(dx)=c_n\int_{U(n)\times U(n)}\int_{\Bbb R_+^n}f(u\lambda v)\,\alpha(du)\alpha(dv)\prod_{j<k}\big(\lambda_j^2-\lambda_k^2\big)^2\prod_{j=1}^n\lambda_j\,d\lambda_j .$$ $c_n$ is the following constant : \\ $$c_n=\frac{2^n\pii^{n^2}}{n!\Big(\prod_{j=1}^{n-1}j!\Big)^2}.$$
\end{pro} 

\bigskip

By using the preceding results, the resolution of the Cauchy problem for the heat equation on $V_n$ leads to the evaluation of the orbital integral $\mathscr{I}(x,y)$, which is defined for $x,y\in V_n$ by  $$\mathscr{I}(x,y)=\int_{U(n)}\int_{U(n)}e^{\mathfrak{Re}\,{\rm tr}(xuyv^*)}\alpha(du)\alpha(dv).$$

\bigskip

One can remark that the function $\mathscr{I}(x,y)$ is determined by its restriction to the subspace of diagonal matrices because it is invariant under $K_n$ : $$\mathscr{I}(uxv^*,y)=\mathscr{I}(x,uyv^*)=\mathscr{I}(x,y)\ \ (u,v\in U(n)).$$
The Cauchy problem for the heat equation $$\dif{U}{t}=\Delta U,$$ $$U(0,x)=f(x),$$ where $f$ is a bounded continuous function on $V_n$, has a unique solution which is $$ U(t,x)=\frac{1}{(4\pi t)^{\frac{N}{2}}}\int_{V_n}e^{-\frac{1}{4t}|||x-y|||^2}f(y)\,m(dy) \ \ (t>0,\, x\in V_n),$$ where $N=2n^2$ is the dimension of $V_n$, $|||.|||$ is the Hilbert-Schmidt norm on $V_n$ and $m$ is the Euclidean measure.\\

Let us assume that the function $f$ is invariant under the action of $K_n$. Then the solution $U$ will be also $K_n$-invariant. Hence we can write $$f(x)=f_0(\lambda),\ \ U(t,x)=U_0(t,\lambda).$$ By using the Weyl integration formula (Proposition \ref{intprop}), the solution $U_0(t,\lambda)$ is given by:

$$U_0(t,\lambda)=\int_{\Bbb R_+^n} H_0(t,\lambda,\theta)f_0(\theta)D(\theta)\prod_{j=1}^n\theta_j\,d\theta_j,$$ 
with
\ba H_0(t,\lambda,\theta)&=& c_n\frac{1}{(4\pi t)^{\frac{N}{2}}}\int_{U(n)}\int_{U(n)} e^{-\frac{1}{4t}|||\lambda-u\theta v^*|||^2}f(y)\,\alpha(du)\alpha(dv)\\ 
&=& c_n\frac{1}{(4\pi t)^{\frac{N}{2}}}e^{-\frac{1}{4t}(||\lambda||^2+||\theta||^2)}\int_{U(n)}\int_{U(n)}e^{\frac{1}{2t}\mathfrak{Re}\,{\rm tr}(\lambda u\theta v^*)}\alpha(du)\alpha(dv)\\
&=&c_n\frac{1}{(4\pi t)^{\frac{N}{2}}}e^{-\frac{1}{4t}(||\lambda||^2+||\theta||^2)}\mathscr{I}(\frac{1}{2t}\lambda,\theta).\ea

\bigskip

\begin{teo} If $\lambda={\rm diag(}\lambda_1,\dots,\lambda_n{\rm)}$ and $\theta={\rm diag(}\theta_1,\dots,\theta_n{\rm)},$ then $$\mathscr{I}(\lambda,\theta) =2^{{n(n-1)}}\big[1!2!\times\dots\times(n-1)!\big]^2\,\frac{1}{D(\lambda)D(\theta)}\,{\rm det\Big(}\big(I_0(\lambda_i\theta_j)\big)_{1\leq i,j\leq n}{\rm \Big)},$$ where $I_0$ is the modified Bessel function: $$I_0(z)=\sum_{k=0}^\infty\frac{1}{2^{2k}(k!)^2}\,z^{2k} \ \ (z\in\Bbb C).$$
\end{teo}
\begin{paragraph}{\bfseries{Proof.}}{\rm From the formula that gives the radial part of the laplacian (Proposition \ref{radprop}) one can deduce that the function $U_0$ is a solution of the equation\\ $$\diff{U_0}{t}=\frac{1}{D(\lambda)}\,\sum_{i=1}^n\Big(\diff{}{\lambda_i}+\frac{1}{\lambda_i}\dif{}{\lambda_i}\Big)\big(D(\lambda)F(\lambda)\big).$$ We put then $$V(t,\lambda)=D(\lambda)U_0(t,\lambda),\ \ g(\lambda)=D(\lambda)f_0(\lambda).$$

\bigskip

 The function $V$ is in consequence a solution of the Cauchy problem\\
 
  $$\diff{V}{t}=\sum_{i=1}^n\Big(\diff{V}{\lambda_i}+\frac{1}{\lambda_i}\dif{V}{\lambda_i}\Big),$$ $$V(0,\lambda)=g(\lambda).$$ 
  
  \bigskip

Let us assume that $f$ belongs to the Schwartz space $\mathscr{S}(V_n)$. By a result similar to (\cite{Far4}, Lemme X-3.1) one can prove that, for every $T>0$, the function $V$ is bounded on $[0,T]\times\Bbb R_+^n$.

 On the other hand, for $n=1$, the preceding problem is equivalent, in cylindrical coordinates, to the Cauchy problem for the heat equation on $\Bbb R^2$ with $f$ radial. The solution of such a problem is expressed using the modified Bessel function $I_0$ and it is given by (see \cite{Far4}, Chapter IX, exercice 3) :
 
  $$\frac{1}{2t}\int_0^\infty e^{-\frac{r^2+\rho^2}{4t}}g(\rho)\,I_0\Big(\frac{r\rho}{2t}\Big)\,\rho\,d\rho.$$ 
  
  \bigskip
  
  In consequence we can conclude that the solution of our problem for an arbitrary $n$ is given by $$V(t,\lambda)=\frac{1}{(2t)^n}\int_{\Bbb R_+^n}e^{-\frac{1}{4t}(||\lambda||^2+||\theta||^2)}g(\theta)\prod_{i=1}^nI_0\Big(\frac{\lambda_i\theta_i}{2t}\Big)\,\theta_i\,d\theta_i.$$  As the function $g$ is skewsymmetric this last relation can be written $$V(t,\lambda)=\frac{1}{(2t)^n}\int_{\Bbb R_+^n}e^{-\frac{1}{4t}(||\lambda||^2+||\theta||^2)}g(\theta)\frac{1}{n!}\sum_{\sigma\in\mathfrak{S}_n}\varepsilon(\sigma)\prod_{i=1}^nI_0\Big(\frac{\lambda_i\theta_{\sigma(i)}}{2t}\Big)\,\theta_i\,d\theta_i.$$ Hence, for every function $g(\theta)=D(\theta)f_0(\theta)$, where $f_0$ is a symmetric function in $\mathscr{S}(\Bbb R_+^n),$ \ba\lefteqn{\int_{\Bbb R_+^n} H_0(t,\lambda,\theta)g(\theta)\prod_{i=1}^n\theta_i\,d\theta_i}\\& & = \frac{1}{n!(2t)^n}\int_{\Bbb R_+^n}e^{-\frac{1}{4t}(||\lambda||^2+||\theta||^2)}g(\theta)\,{\rm det\bigg(}\Big(I_0\big(\frac{\lambda_i\theta_j}{2t}\big)\Big)_{1\leq i,j\leq n}{\rm \bigg)}\,\prod_{i=1}^n\theta_i\,d\theta_i .\ea This proves that the kernel $H_0$ is equal to $$H_0(t,\lambda,\theta)=\frac{1}{n!(2t)^n}\,\frac{1}{D(\lambda)D(\theta)}\,e^{-\frac{1}{4t}(||\lambda||^2+||\theta||^2)}\,{\rm det\bigg(}\Big(I_0\big(\frac{\lambda_i\theta_j}{2t}\big)\Big)_{1\leq i,j\leq n}{\rm \bigg)}.$$ We obtain the result by comparing, for $t=\frac{1}{2}$, the two expressions we obtained for $H_0$.$\hspace{1cm}\Box$}
\end{paragraph}

\bigskip

For $x\in V_n$ the orbital measure $\m_x$ is defined on $V_n$ by\\ $$\int_{V_n}f(y)\,\m_x(dy)=\int_{U(n)}\int_{U(n)}f(uxv^*)\,\alpha(du)\alpha(dv),$$ where $\alpha$ is the normalized Haar measure on $U(n)$ and $f$ is a continuous function on $V_n$. The Fourier transform of $\m_x$ is the following function $\widehat{\m_{x}}$ :\\ \ba \widehat{\m_x}(\xi)&=&\int_{V_n}e^{i\<\xi,y\>}\m_x(dy)\\
&=&\int_{U(n)}\int_{U(n)}e^{i\mathfrak{Re}\,{\rm tr}(\xi uxv^*)}\alpha(du)\alpha(dv)\\
&=& \mathscr{I}(x,i\xi).\ea   

\bigskip

The spherical functions of positive type for the Gelfand pair $(G_n,K_n)$ are the functions $\pp_x=\widehat{\m_{x}}, \ (x\in\Bbb R^n)$, Fourier transforms of the orbital measures $\m_x$.\\

\begin{cor}\label{spherfin} If $x={\rm diag(}x_1,\dots,x_n{\rm)}$ and $\xi={\rm diag(}\xi_1,\dots,\xi_n{\rm)},$ $$\pp_x(\xi):=\widehat{\m_x}(\xi)=\big({\bf \delta!}\big)^2\,\frac{(-4)^{\frac{n(n-1)}{2}}}{D(x)D(\xi)}\,{\rm det\Big(}\big(J_0(x_j\xi_k)\big)_{1\leq j,k\leq n}{\rm \Big)},$$ where  $${\bf \delta}=(\delta_1,\delta_2,\dots,\delta_n):=(n-1,n-2,\dots,0),$$ $${\bf \delta!}=\delta_1!\times\dots\times\delta_n!,$$ and $J_0$ is the classical Bessel function.
\end{cor}

\bigskip
\section{Multiplicativity property of spherical functions} 

\bigskip

A spherical function, for the spherical pair $(G_\infty,K_\infty)$, is a continuous function $\pp$ on $G_\infty$ satisfying $$\lim_{n\rightarrow\infty}\int_{K_n}\pp(xky)dk=\pp(x)\pp(y),$$ where $dk$ is the normalized Haar measure of the product group $K_n=U(n)\times U(n)$ (see \cite{Far2}, Theorem 5.1). In our case, the function $\pp$ can be seen as a function on $V_\infty$ and hence  $$\lim_{n\rightarrow\infty}\int_{U(n)\times U(n)}\pp(x+k_1yk_2^*)\alpha_n(dk_1)\alpha_n(dk_2)=\pp(x)\pp(y),$$ where $\alpha_n$ is the normalized Haar measure of the unitary group $U(n)$.\\

\begin{teo}\label{multi}{{\rm (The multiplicativity property)}} Let $\pp\in\mathfrak{P}$. The function $\pp$ is spherical, if and only if, there exists a continuous function $\Phi$ on $\Bbb R$, with $\Phi(0)=1$ such that $$\pp\big({\rm diag}(a_1,\dots,a_n,0,\dots)\big)=\Phi(a_1)\dots\Phi(a_n).$$
\end{teo}

\bigskip

Let us put, for $m\leq n$,$$K_m(n)=\Biggl\{\left(\begin{array}{cc}       I_m & 0  \\       0  & v \end{array}\right)\,\Big|\, v\in U(n-m)\Biggl\}\simeq U(n-m), $$ and $$K_m(\infty)=\bigcup_{n=1}^\infty K_m(n)\subset U(\infty).$$ Also, let us put $$K(m,n)=\Biggl\{\left(\begin{array}{cc}       u & 0  \\       0  & v \end{array}\right)\,\Big|\, u\in U(m),\, v\in U(n-m)\Biggl\}\subset U(n).$$ We introduce in addition, for $n\geq 2m$, a Cartan subgroup for the symmetric pair $\big(U(n),\,K(m,n)\big)$:$$a(\theta)=\left(\begin{array}{ccccccc}  \cos\theta_1 & & & -\sin\theta_1 & & & \\
 & \ddots & & & \ddots & & \\
 & & \cos\theta_m & & & -\sin\theta_m & \\
\sin\theta_1 & & & \cos\theta_1 & & & \\
 & \ddots & & & \ddots & & \\
 & & \sin\theta_m & & & \cos\theta_m & \\
 & & & & & & I_{n-2m}\end{array}\right).$$ Hence, for every $k\in K(n)$, $$k=h_1a(\theta)h_2,$$ with $h=(h_1,h_2)\in K(m,n)\times K(m,n)=:K^2(m,n)$, and the Weyl integration formula that corresponds to this last decomposition is given, for every integrable function $f$ on $U(n)\times U(n)$, by  \ba\lefteqn{\int_{U(n)\times U(n)} f(k_1,\,k_2)\,\alpha_n(dk_1)\alpha_n(dk_2)}\\& & =\int_{[0,\pi]^m\times [0,\pi]^m}\int_{K^2(m,n)}\int_{K^2(m,n)}f\big(h_1a(\theta)h_2,\,g_1a(\zeta)g_2\big)\,\kappa(dh)\kappa(dg)D_{m,n}(\theta)D_{m,n}(\zeta)\,d\theta\,d\zeta,\ea where $\kappa=\beta\otimes\beta$ and $\beta$ is the standardized Haar measure of $K(m,n)$, and $$D_{m,n}(\theta)=c_{m,n}\Big|\prod_{1\leq i<j\leq m}\big(\sin(\theta_i+\theta_j)\big)^2\big(\sin(\theta_i-\theta_j)\big)^2\prod_{i=1}^m\big(\sin2\theta_i\big)\big(\sin\theta_i\big)^{2(n-2m)}\Big|,$$ where $c_{m,n}$ is a constant such that $$\int_{[0,\pi]^m}D_{m,n}(\theta)\,d\theta_1\dots\theta_m=1.$$

\bigskip

 \begin{pro}\label{convkey} Let $f$ be a continuous function on $K_\infty$ which is $K_m(\infty)\times K_m(\infty)$-invariant. Then \ba\lefteqn{\lim_{n\rightarrow\infty}\int_{U(n)\times U(n)}f(k_1,\,k_2)\,\alpha_n(dk_1)\alpha_n(dk_2)}\\& & =\int_{K_m^2}\int_{K_m^2}f(h_1w_mh_2,\,g_1w_mg_2)\,\alpha_m(dh_1)\alpha_m(dh_2)\alpha_m(dg_1)\alpha_m(dg_2),\ea
  where $$w_m=a(\frac{\pi}{2},\dots,\frac{\pi}{2})=\left(\begin{array}{ccc}  0 & -I_m & 0 \\
I_m &  0   & 0 \\
 0 & 0 & I \end{array}\right).$$
\end{pro}

\bigskip

\begin{lem}\label{convmult} Let $X$ be a compact topological space and $\m$ a positive measure on $X$ such that for every non-null open set $B$ we have $\m(B)\geq0$. Let $\delta$ be a continuous positive function on $X$, which reaches its maximum at a unique point $x_0$. Let us put $$\frac{1}{c_n}=\int_X\delta(x)^n\,\m(dx).$$ Then, if $f$ is a continuous function on $X\times X$, $$\lim_{n\rightarrow\infty}c_n^2\int_{X\times X}f(x,y)\delta(x)^n\delta(y)^n\,\m(dx)\m(dy)=f(x_0,x_0).$$
\end{lem}
\begin{paragraph}{\bfseries{Proof.}}{\rm We apply Lemma 5.4 in \cite{Far2} to the function $\delta(x,y)=\delta(x)\delta(y)$.$\hspace{1cm}\Box$}
\end{paragraph}

\bigskip

By using Lemma \ref{convmult}, for every fixed $m$, and for every continuous function $f$ on $[0,\pi]^m\times[0,\pi]^m$,$$\lim_{n\rightarrow\infty}\int_{[0,\pi]^m\times[0,\pi]^m}f(\theta,\,\zeta)D_{m,n}(\theta)D_{m,n}(\zeta)\,d\theta \,d\zeta=f\Big(\frac{\pi}{2},\dots,\frac{\pi}{2};\frac{\pi}{2},\dots,\frac{\pi}{2}\Big).$$

\bigskip

\begin{paragraph}{\bfseries{Proof of Proposition \ref{convkey}.\\}}{\rm By using the integration formula and the invariance under $K_m(\infty)\times K_m(\infty)$, we obtain $$\int_{U(n)\times U(n)} f(k_1,\,k_2)\,\alpha_n(dk_1)\alpha_n(dk_2)=\int_{[0,\pi]^m\times [0,\pi]^m}F(\theta,\,\zeta)D_{m,n}(\theta)D_{m,n}(\zeta)\,d\theta\,d\zeta,$$ with $$F(\theta,\,\zeta)=\int_{K_m^2}\int_{K_m^2}f\big(h_1a(\theta)h_2,\,g_1a(\zeta)g_2\big)\,\alpha_m(dh_1)\alpha_m(dh_2)\alpha_m(dg_1)\alpha_m(dg_2).$$ As a result, by Lemma \ref{convmult}  \ba\lefteqn{\lim_{n\rightarrow\infty}\int_{U(n)\times U(n)} f(k_1,\,k_2)\,\alpha_n(dk_1)\alpha_n(dk_2)}\\& &=\int_{K_m^2}\int_{K_m^2}f(h_1w_mh_2,\,g_1w_mg_2)\,\alpha_m(dh_1)\alpha_m(dh_2)\alpha_m(dg_1)\alpha_m(dg_2).\hspace{1cm}\Box\ea}

\end{paragraph}

\bigskip

\begin{cor}\label{cormult} Let $\pp$ be a $K_\infty$-invariant continuous function on $V_\infty$. Then, for\\ $x={\rm diag(}a_1,\dots,a_m,0,\dots{\rm )}$ and $y={\rm diag(}b_1,\dots,b_m,0,\dots{\rm )}$, $$\lim_{n\rightarrow\infty}\int_{U(n)\times U(n)}\pp(x+k_1yk_2^*)\,\alpha_n(dk_1)\alpha_n(dk_2)=\pp\big({\rm diag(}a_1,\dots,a_m,b_1,\dots,b_m,0\dots{\rm )}\big).$$
\end{cor}
\begin{paragraph}{\bfseries{Proof.}}{\rm The function $(k_1,\,k_2)\mapsto\pp(x+k_1yk_2^*)$ is $K_m(\infty)\times K_m(\infty)$-invariant. Hence, we can apply Proposition \ref{convkey}: \ba\lefteqn{\lim_{n\rightarrow\infty}\int_{U(n)\times U(n)}\pp(x+k_1yk_2^*)\,\alpha_n(dk_1)\alpha_n(dk_2)}\\ & &=\int_{K_m^2}\int_{K_m^2}\pp(x+h_1w_mh_2yg_2^*w_m^{-1}g_1^*)\,\alpha_m(dh_1)\alpha_m(dh_2)\alpha_m(dg_1)\alpha_m(dg_2).\ea Finally, we obtain the result by using the fact that $$x+h_1w_mh_2yg_2^*w_m^{-1}g_1^*\in U(\infty){\rm diag(}a_1,\dots,a_m,b_1,\dots,b_m{\rm )}U(\infty).\hspace{1cm}\Box$$}
\end{paragraph}

\bigskip

\begin{paragraph}{\bfseries{Proof of Theorem \ref{multi}.}}{\rm Let $\pp\in\mathfrak{P}$. If $\pp$ is spherical then for $$x={\rm diag(}a_1,\dots,a_m,0,\dots{\rm )},\  y={\rm diag(}b_1,\dots,b_m,0,\dots{\rm )},$$ $$\lim_{n\rightarrow\infty}\int_{U(n)\times U(n)}\pp(x+k_1yk_2^*)\,\alpha_n(dk_1)\alpha_n(dk_2)=\pp(x)\pp(y).$$ By Corollary \ref{cormult}, $$\pp\big({\rm diag(}a_1,\dots,a_m,0,\dots{\rm )}\big)\big({\rm diag(}b_1,\dots,b_m,0,\dots{\rm )}\big)=\pp\big({\rm diag(}a_1,\dots,a_m,b_1,\dots,b_m,0,\dots{\rm )}\big).$$ By applying Corollary \ref{cormult} many times as necessary, one obtains $$\pp\big({\rm diag(}a_1,\dots,a_n,0,\dots{\rm )}\big)=\Phi(a_1)\dots\Phi(a_n),$$ where $$\Phi(\lambda)=\pp\big({\rm diag(}\lambda,0,\dots{\rm )}\big).$$

Inversely, let us assume that there exists a continuous function $\Phi$ on $\Bbb R$ such that $$\pp\big({\rm diag(}a_1,\dots,a_n,0,\dots{\rm )}\big)=\Phi(a_1)\dots\Phi(a_n).$$ Then, by Corollary \ref{cormult}, $$\lim_{n\rightarrow\infty}\int_{U(n)\times U(n)}\pp(x+k_1yk_2^*)\,\alpha_n(dk_1)\alpha_n(dk_2)=\pp(x)\pp(y).$$ Hence the function $\pp$ is spherical.$\hspace{1cm}\Box$ \\}
\end{paragraph}

\bigskip

\bigskip


\section{Modified P{\'o}lya functions : definition and convergence}

\bigskip

\begin{defi}\label{dif}\rm{
{\it The modified P{\'o}lya function} of parameter $\omega=(\alpha,\gamma)$ with $\alpha=\{\alpha_j\}_{j\geq1}\in\ell^{1}(\Bbb N)$, $\alpha_j\in\Bbb R_+$ and $\gamma\in\Bbb R_+$ is defined on $\Bbb R$ by :\\
$${\Pi}(\omega,\lambda):=e^{-\frac{1}{4}\gamma\lambda^2}\prod_{j=1}^{\infty}\frac{1}{1+\frac{1}{4}\alpha_j
  \lambda^2}.$$}
\end{defi}

\bigskip

We consider on the set $\mathcal P$ of modified P{\'o}lya functions the topology of uniform convergence on compact sets of $\Bbb R$. The topological space $\mathcal P$ is metrizable and complete. This topology can be expressed in terms of the set of parameters : $$\Omega=\left\{\omega=(\alpha,\gamma)\,\Big|\, \alpha=(\alpha_j)_{j\geq1},\, \alpha_j\geq0,\, \sum_{j=1}^\infty\alpha_j<\infty,\, \gamma\geq0\right\}.$$ For a continuous function $f$ on $\Bbb R$, we define the function $L_f$ on $\Omega$ by \bq\label{log}L_f(\omega)=\int_{\mathbb R} f(t)\sigma_\omega(dt)=\gamma f(0)+\sum\limits_{j=1}^{\infty}
{\alpha}_j f(\alpha_j).\eq Let us remark that the moments of the measure $\sigma_\omega$ are given by $${\mathcal M}_0(\sigma_\omega)=\int_{\Bbb
  R}\sigma_\omega(dt)=\gamma+\sum_{k=1}^{\infty}\alpha_k=\gamma+p_1(\alpha),$$
and for $m\geq 1$,
$${\mathcal M}_m(\sigma_\omega)=\int_{\Bbb
  R}t^m\sigma_\omega(dt)=\sum_{k=1}^{\infty}\alpha_k^{m+1}=p_{m+1}(\alpha),$$
where $p_m$ is the Newton power sum function : for $x=(x_1,x_2,\ldots)\in\ell^{1}(\Bbb N)$ and $m\geq 1$, $$p_m(x)=\sum_{k=1}^{\infty}x_k^m.$$

\bigskip

We consider on $\Omega$ the initial topology associated to the functions $L_f$. A point $\omega\in \Omega$ is seen as a point configuration, i.e. a permutation of the numbers $\{\alpha_k\}$, $\gamma$ does not change $\omega$. For $\lambda$ fixed, the function $\omega\mapsto \Pi(\omega,\lambda)$ is injective and continuous on $\Omega$. This can be seen by looking at the logarithmic derivative of $\Pi(\omega,\lambda)$ : \bq\label{logg}\begin{aligned}{{\Pi^{'}(\omega,\lambda)}\over
  {\Pi(\omega,\lambda)}}&=-\dfrac{1}{2}(\gamma+p_{1}(\alpha))\lambda+i\sum\limits_{m=2}^{\infty}p_m(\alpha)\left(\dfrac{i\lambda}{2}\right)^{2m-1}.\end{aligned}\eq 

\bigskip

\begin{lem}\label{reff} Let $\mathcal I$ be the set of positive measures $\m$ on $\Bbb R_+^*$ such that $\m\left([a,+\infty[\right)\in\Bbb N$, for all $a>0$. Then, the following properties hold :\\

 {\rm (i)} For all $\m\in\mathcal I$, there exists a sequence of positive reals $\{\alpha_k\}_k$ such that $$\m=\sum_{k=1}^N\delta_{\alpha_k},$$ with $N\leq\infty$. If $N=\infty$, then the sequence $\{\alpha_k\}$ converges to $0$.\\

 {\rm (ii)} Let $\m_n$ be a sequence of measures in $\mathcal I$. Assume that there exists a measure $\m$ on $\Bbb R_+^*$ such that, for every function $f$ in the set $\mathcal C_b^0(\Bbb R_+^*)$ of bounded continuous functions on $\Bbb R_+^*$, vanishing near $0$, \bq\label{us}\lim_{n\to\infty}\int_{\Bbb R_+^*}f(x)\;\m_n(dx)=\int_{\Bbb R_+^*}f(x)\;\m(dx).\eq Then $\m\in\mathcal I$.\\

 \end{lem}
 
 \begin{paragraph}{\bfseries{Proof.}}{\rm (i) For $t>0$, let $F$ be the function defined on $\Bbb R_+^*$ by $F(t)=\m\left([t,\infty[\right)$. We can observe that $F$ is integer-valued, decreasing, and left-continuous. In consequence, the set $D$ of discontinuity points of $F$ is countable : $D=\left\{\alpha_k,\ k\in I\right\}$. Furthermore, we can observe that $D=\{\ t>0 \ |\ \m\left(\{t\}\right)\neq0\}$. We can also remark that the jump at each discontinuity point $\alpha_k$ is an integer $m_k$. Hence, if the number of discontinuity points $N$ is finite, then the measure $\m$ is given by : $\m=\sum_{k=1}^N m_k\delta_{\alpha_k}$ and, provided that we repeat the $\alpha_k$ as many times as their multiplicity $m_k$, we get $\m=\sum_{k=1}^N\delta_{\alpha_k}$. Finally, if $D$ is not finite, i.e. $N=\infty$, then the sequence $\{\alpha_k\}$ converges necessarily to $0$. 
 
 (ii) Let $F_n$ be the function defined on $\Bbb R_+^*$ by $F_n(t)=\m_n\left([t,\infty[\right)$, $F(t)=\m\left([t,\infty[\right)$, and $D$ the set of discontinuity points of $F$. For $t\in D^c$ fixed, the positive measures $\m_n$ and $\m$ are bounded on $[t,\infty[$. In addition, we have $\m(\{t\})=\m(\{\infty\})=0$. Furthermore, it follows by (\ref{us}) that $\m_n$ converges weakly to $\m$ on $[t,\infty[$. As a result $$\lim_{n\to\infty}F_n(t)=F(t)\quad (t\in D^c).$$ Since $F_n(t)$ is a sequence of integers, it follows that, if $t\in D^c$, then $F(t)\in\Bbb N$. On the other hand, because $D$ is countable, its complementary set $D^c$ is dense in $\Bbb R_+^*$. This implies, since $F$ is left-continuous, that $F(t)\in\Bbb N$ for all $t>0$. Hence $\m\in\mathcal I$. $\hspace{1cm}\Box$\\ }

\end{paragraph}

\bigskip

\begin{teo}\label{fermx} $\Omega$ is weakly closed in the set of bounded measures on $\Bbb R_+$.
\end{teo}

\begin{paragraph}{\bfseries{Proof.}}{\rm Let $\omega^{(n)}=(\alpha^{(n)},\gamma^{(n)})$ be a sequence in $\Omega$. Assume that there is a bounded positive measure $\sigma$ on $\Bbb R_+$ such that, for every bounded continuous function $f$ on $\Bbb R_+$, \bq\label{cutis}\lim_{n\to\infty}\gamma^{(n)}f(0)+\sum_{k=1}^\infty\alpha_k^{(n)}f(\alpha_k^{(n)})=\int_{\Bbb R_+}f(t)\sigma(dt) .\eq For every $a>0$, there is a finite number of $\alpha_k^{(n)}$ which are greater than $a$. As a result, the measure $\m_n$ on $\Bbb R_+^*$ given by $$\int_{\Bbb R_+^*}f(t)\;\m_n(dt)=\sum_{k=1}^\infty f(\alpha_k^{(n)}) $$ belongs to $\mathcal I$. Furthermore, by (\ref{cutis}), we can observe that $\m_n$ converges weakly on $\Bbb R_+^*$ to $\dfrac{\tilde{\sigma}}{t}$, where $\tilde{\sigma}$ is the restriction of the measure $\sigma$ to $\Bbb R_+^*$. By (ii) of Lemma \ref{reff}, there exists a sequence $\alpha=\{\alpha_k\}_k$ such that, for every function $f\in\mathcal C_b^0(\Bbb R_+^*)$, $$\lim_{n\to\infty}\int_{\Bbb R_+^*}f(t)\;\m_n(dt)=\sum_{k=1}^\infty f(\alpha_k).$$ Since thet set $\mathcal C_c(\Bbb R_+^*)$ of continuous functions with compact support in $\Bbb R_+^*$ is included in $\mathcal C_b^0(\Bbb R_+^*)$, it follows that $\tilde{\sigma}$ is given by $$\int_{\Bbb R_+^*}f(t)\;\tilde{\sigma}(dt)=\sum_{k=1}^\infty\alpha_kf(\alpha_k).$$ Since the measure $\sigma$ is bounded, $\sum_{k=1}^\infty\alpha_k<\infty$. Finally $$\sigma=\gamma\delta_0+\sum_{k=1}^\infty\alpha_k\delta_{\alpha_k},$$ with $\gamma=\sigma(\{0\})$. Hence $\sigma=\sigma_\omega$ with $\omega=(\alpha,\gamma)\in\Omega$. $\hspace{1cm}\Box$ \\}

\end{paragraph}

\bigskip

A modified P{\'o}lya function of parameter $\omega$ is the Fourier transform of a probability measure $\m_\omega$. Let $\mathfrak M_\Omega$ be the set of these measures : $$\mathfrak M_\Omega=\left\{\m_\omega \ |\ \Pi(\omega,\lambda)=\widehat{\m_\omega}(\lambda)\right\}.$$ We consider on $\mathfrak M_\Omega$ the weak topology of measures. We will prove that the topology defined on $\Omega$ is equivalent to the weak topology of $\mathfrak M_\Omega$. We will need the following proposition (see \cite{Far8}, Proposition 3.11) :\\

\begin{pro}\label{proofc} Let $\psi_n$ be a sequence of $\mathcal C^\infty$-functions on $\Bbb R^d$ of positive type with $\psi_n(0)=1$, and $\psi$ an analytic function on a neighborhood of $0$. Assume that, for every $\alpha=(\alpha_1,\dots,\alpha_d)$, $$\lim_{n\to\infty}\partial^\alpha\psi_n(0)=\partial^\alpha\psi(0).$$ Then $\psi$ has an analytic extension to $\Bbb R^d$, and $\psi_n$ converges to $\psi$ uniformly on compact sets in $\Bbb R^d$. \\\end{pro}

\bigskip

\begin{pro}\label{irman} The topology of $\Omega$ is equivalent to the topology of $\mathfrak M_\Omega$.\\
\end{pro}

\begin{paragraph}{\bfseries{Proof.}}\rm{ (i) Assume that $\omega^{(n)}$ converges to $\omega$ in the topology of $\Omega$. Then the modified P{\'o}lya functions $\Pi(\omega^{(n)},\lambda)$ and $\Pi(\omega,\lambda)$ are holomorphic in $D(0,R)$ where $\displaystyle\frac{1}{R}=\sup\limits_{m,
  n}\left|\alpha_m^{(n)}\right|$. Therefore, for every $\lambda\in D(0,R)$, their logarithmic derivatives are given by $${{\Pi^{'}(\omega^{(n)},\lambda)}\over {\Pi(\omega^{(n)},\lambda)}}=-\dfrac{1}{2}(\gamma^{(n)}+p_1(\alpha^{(n)}))\lambda+i\sum\limits_{m=2}^{\infty}p_m(\alpha^{(n)})\left(\dfrac{i\lambda}{2}\right)^{2m-1},$$
$${{\Pi^{'}(\omega,\lambda)}\over
  {\Pi(\omega,\lambda)}}=-\dfrac{1}{2}(\gamma+p_1(\alpha))\lambda+i\sum\limits_{m=2}^{\infty}p_m(\alpha)\left(\dfrac{i\lambda}{2}\right)^{2m-1}.$$
For every bounded continuous function $f$ on $\Bbb R$ $$\lim_{n\to\infty}\int_{\Bbb R} f(t)\sigma_{\omega^{(n)}}(dt)=\int_{\Bbb R}
f(t)\sigma_\omega(dt),$$
where $\sigma_{\omega^{(n)}}$ is the bounded positive measure on $\Bbb R$
  associated to $\omega^{(n)}=(\alpha^{(n)},\gamma^{(n)})$.
Hence, for $f=1$, the sequence $\gamma^{(n)}+p_1(\alpha^{(n)})$ converges. In consequence, there exists a constant $A>0$ such that, for every $n$,
$$ 0\leq \gamma^{(n)}+p_1(\alpha^{(n)})\leq
A.$$\\
Therefore, the sequence $\{\alpha_m^{(n)}\}$ is bounded by $A$ for every $m$ and every $n$.\\ Hence,
$$\mbox{supp}(\sigma_{\omega^{(n)}})\subset [-A,A]\,\mbox{ and}\ \
\mbox{supp}(\sigma_\omega)\subset [-A,A].$$
It follows that the sequence of measures $\sigma_{\omega^{(n)}}$ converges to $\sigma_\omega$ for every continuous function on $\Bbb R$. We can then deduce that, for every $m\geq
2$,$$\begin{aligned}\lim\limits_{n\to {\infty}}p_m({{\alpha}^{(n)}})&=\lim\limits_{n\to
  {\infty}}\int_{\mathbb R}t^{m-1}\sigma_{\omega^{(n)}}(dt)\\&=\int_{\mathbb
  R}t^{m-1}\sigma_\omega(dt)=p_m(\alpha),\end{aligned}$$ and $$p_{m}({{\alpha}^{(n)}})\leq
\big(p_{1}{({\alpha}^{(n)})}\big)^{m}\leq A^{m}.$$\\
Furthermore, we have
$$\lim_{n\to\infty}\gamma^{(n)}+p_1(\alpha^{(n)})= \gamma+p_1(\alpha).$$
 In consequence,
  $$\lim_{n\to\infty}{{\Pi^{'}(\omega^{(n)},\lambda)}\over {\Pi(\omega^{(n)},\lambda)}}={{\Pi^{'}(\omega,\lambda)}\over {\Pi(\omega,\lambda)}},$$
 the convergence being uniform on compact sets in $D(0,R)$. Hence, for $|\lambda|< R$,
  $$\lim_{n\to\infty}\Pi(\omega^{(n)},\lambda)=\Pi(\omega,\lambda),$$ since
  $\Pi(\omega^{(n)},0)=1$ and $\Pi(\omega,0)=1$.\\

 The functions
  $\Pi(\omega^{(n)},\lambda)$ and $\Pi(\omega,\lambda)$ being of positive type, by Proposition \ref{proofc}, $\Pi(\omega^{(n)},\lambda)$ converges uniformly on compact sets in $\Bbb R$ to $\Pi(\omega,\lambda)$. Finally, by applying the L\'evy-Cramer theorem, one can prove that $\m_{\omega^{(n)}}$ converges weakly to $\m_\omega$.\\

(ii) Assume that $\m_{\omega^{(n)}}$ converges weakly to $\m_\omega$. This implies that $\Pi(\omega^{(n)},\lambda)$ converges uniformly on compact sets in $\Bbb R$ to $\Pi(\omega,\lambda)$. Let $\lambda_0>0$. Since the modified P{\'o}lya function $\Pi(\omega^{(n)},\lambda)$ is continuous, non-zero on $\Bbb R$ and satisfies $\Pi(\omega^{(n)},0)=1$, there exists $M>0$ such that, for every $n$, $\Pi(\omega^{(n)},\lambda_0)\geq M.$ In consequence $$\dfrac{1}{4}\gamma^{(n)}\lambda_0^2+\dfrac{1}{4}\lambda_0^2\sum_{k=1}^\infty\alpha_k^{(n)}\leq e^{\frac{1}{4}\gamma^{(n)}\lambda_0^2}\prod_{k=1}^\infty\left(1+\dfrac{1}{4}\alpha_k^{(n)}\lambda_0^2\right)\leq {1\over M}.$$ Therefore $$-2\Pi^{''}(\omega^{(n)},0)=p_1(\alpha^{(n)})+\gamma^{(n)}\leq{4\over{\lambda_0^2\,M}}=:R.$$ 

\bigskip 

It follows that the function $\Pi(\omega^{(n)},\lambda)$ is holomorphic for $|\lambda|<R$ and then also in the strip $\Sigma_R=\{x+iy\,|\,|y|<R\}$ (see \cite{Far2}, Lemma 6.5). Furthermore, for $r<R$, there exists a constant $M(r)>0$ such that $|\Pi(\omega^{(n)},\lambda)|\leq M(r) \ {\rm for}\ \lambda\in\Sigma_r.$ From the theorem of Montel, it follows that there is a subsequence $\Pi(\omega^{(n_j)},\lambda)$ which converges uniformly on compact sets in $\Sigma_R$. Since the sequence itself converges to $\Pi(\omega,\lambda)$ on $\Bbb R$, we get that $\Pi(\omega^{(n)},\lambda)$ converges to $\Pi(\omega,\lambda)$ uniformly on compact sets in $\Sigma_R$. As a result, the logarithmic derivatives $\frac{\Pi^{'}(\omega^{(n)},\lambda)}{\Pi(\omega^{(n)},\lambda)}$ converge uniformly in a neighborhood of $0$ which implies the convergence of the coefficients of their Taylor expansions at $0$ : $$\lim_{n\to\infty}\gamma^{(n)}+p_1(\alpha^{(n)})=\gamma+p_1(\alpha),$$ $$\lim_{n\to\infty}p_m(\alpha^{(n)})=p_m(\alpha), \ \ {\rm for} \ m\geq2.$$ This proves that $\omega^{(n)}$ converges to $\omega$ in the topology of $\Omega$. $\hspace{1cm} \Box$\\ }
\end{paragraph}

\bigskip

Let us put for every modified P{\'o}lya function of parameter $\omega=(\alpha,\gamma)$,\\ $$\widetilde{p}_1(\omega)=p_1(\alpha)+\gamma,$$ and, for every $R\geq0$,\\ $$\Omega_R=\{\omega\in\Omega \ |-2\Pi^{''}(\omega,0)=\widetilde{p}_1(\omega)\leq R\}.$$

\bigskip

\begin{cor}\label{compact} $\Omega_R$ is a compact subset of $\Omega$. \end{cor}
\begin{paragraph}{\bfseries{Proof.}}\rm{ Since the modified P{\'o}lya function is the Fourier transform of a probability measure $\m$, we have  $$-\Pi^{''}(\omega,0)=\int_{\Bbb R}t^2\m_\omega(dt).$$ We also know that the set of probability measures such that $\int_{\Bbb R}t^2\m_\omega(dt)\leq {R\over2}$ is relatively compact. By using Theorem \ref{fermx} and Proposition \ref{irman}, we can conclude that the set $\Omega_R$ is relatively compact in $\Omega$. Moreover, the convergence of a sequence $\omega^{(n)}$ to $\omega$ in $\Omega$ implies that $\widetilde{p}_1(\omega^{(n)})$ converges to $\widetilde{p}_1(\omega)$. Hence, the map $\omega\mapsto \widetilde{p}_1(\omega)$ is continuous. In consequence, the set $\Omega_R$ is closed and relatively compact in $\Omega$, therefore compact. $\hspace{1cm} \Box$ }
\end{paragraph}

\bigskip

\bigskip
\section{Convergence of orbital measures and spherical functions}  

\bigskip 

Let $(X, {\mathcal B})$ be a measurable space on which a group $G$ acts by measurable transformations. Let $\nu$ be a $G$-invariant probability measure on $X$. A set $E\in\mathcal B$ is said to be $G$-invariant relatively to $\nu$ if, for every $g\in G$, $\nu((gE)\Delta E)=0,$ where $\Delta$ is the symmetric difference. The measure $\nu$ is said to be ergodic relatively to the action of $G$ if, for every $E\in\mathcal
B$ which is $G$-invariant : $\nu(E)=0\ \mbox{or}\ 1.$ If $X$ is locally compact and $K$ a compact group acting on $X$ by homeomorphisms, then the ergodic measures are exactly the orbital ones. It is the case, if one takes $X=V_n=M(n,\Bbb C)$ and $K=K_n=U(n)\times U(n)$.\\

Let $V^{\infty}$ be the space of infinite complex matrices. It is the dual space of $V_\infty$. The space $V^{\infty}$ is defined as the projective limit of the sequence $V_n$ relatively to the orthogonal projections $$p_{m,n}:V_n\rightarrow V_m\ \ (n>m).$$ The space $V^\infty$ can be identified to $\Bbb R^\infty$. The group $K_\infty$ acts on it in the same way as on $V_\infty$. Let $\mathfrak M$ be the set of probability measures on $V^\infty$ which are invariant by $K_\infty$. The Fourier transform of a measure $\nu$ in $\mathfrak M$ is the $K_\infty$-invariant function of positive type $\varphi$ defined on $V_\infty$ by\\
$$\varphi(\xi)=\int_{V^\infty}e^{i\<x,\xi\>}\nu(dx).$$ The Fourier transform establishes a bijection from $\mathfrak{M}$ onto $\mathfrak{P}$ and also between their corresponding subsets of extreme points ext($\mathfrak{M}$) and ext($\mathfrak{P}$). A function $\pp\in\mathfrak{P}$ is spherical if and only if the measure $\nu$ is ergodic.\\

Let $\lambda^{(n)}$ be a sequence of diagonal matrices \\
$$\lambda^{(n)}={\rm diag}({\lambda_1^{(n)}},\ldots,{\lambda_n^{(n)}}).$$ We associate to it the sequence of orbital measures $\nu^{(n)}$ such that, for every continuous function $f$ on $V^\infty$ :\\
$$\int_{V^\infty}f(x)\nu^{(n)}(dx)=\int_{U(n)\times U(n)}f(u\lambda^{(n)}v^*)\alpha_n(du)\alpha_n(dv).$$

\bigskip

The Fourier transform $\pp_n$ of $\nu^{(n)}$ is defined on $V_\infty$ by\\
$$\begin{aligned}\pp_n(\xi)=\pp_n(\xi,\lambda^{(n)})&=\int_{V^\infty}e^{i\<\xi,x\>}\nu^{(n)}(dx)\\
&=\int_{U(n)\times U(n)}e^{i\<\xi, u \lambda^{(n)}v^*\>}\,\alpha_n(du)\alpha_n(dv).\end{aligned}$$

\bigskip

We calculated the preceding integral in Corollary \ref{spherfin}. Let us put $$\ F(z)=\sum_{k=0}^\infty \frac{1}{(k!)^2}\,z^k,\ (z\in\Bbb C),$$ $$\Xi=\big(\Xi_1,\Xi_2,\dots,\Xi_n,0,\dots\big):=\big(-\frac{\xi_1^2}{4},-\frac{\xi_2^2}{4},\dots,-\frac{\xi_n^2}{4},0,\dots\big),$$ $$\Lambda=\big(\Lambda_1,\Lambda_2,\dots,\Lambda_n\big):=\big({\lambda_1^{(n)}}^2,{\lambda_2^{(n)}}^2,\dots,{\lambda_n^{(n)}}^2\big).$$ 

\bigskip

For $\xi={\rm diag}(\xi_1,\dots,\xi_n,0,\dots)$, we can rewrite the spherical function $\pp_n$ as 

\bq\label{formule}\pp_n(\xi,\lambda^{(n)})=\Psi(\Lambda,\ \Xi)=\big({\bf \delta}!\big)^2\,\frac{{\rm det}\bigg(\Big(F\big(\Lambda_i\Xi_j\big)\Big)_{1\leq i,j\leq n}\bigg)}{D(\Lambda)D(\Xi)},\eq where   $$D\big(\Xi\big)=\prod_{i<j}\big(\Xi_i-\Xi_j\big), \ D\big(\Lambda\big)=\prod_{i<j}\big(\Lambda_i-\Lambda_j\big).$$

\bigskip

By a result in (\cite{Far2}, section 2.3), we can write (\ref{formule}) as \bq\label{nes3}\pp_n(\xi,\lambda^{(n)})=\sum_{m_1\geq\dots\geq m_n\geq0} \left(\frac{{\bf \delta!}}{({\bf m+\delta})!}\right)^2\ s_{\bf m}(\Lambda)s_{\bf m}(\Xi).\eq

\bigskip

The generating function of the complete symmetric function $$h_m(x)=\sum_{|\alpha|=m}x^\alpha$$ is given by : $$H(x,t):=\sum_{m=0}^\infty h_m(x)t^m=\prod_{j=1}^n\dfrac{1}{1-x_jt}.$$ The logarithmic derivative of $H(x,-\frac{t^2}{4})$ with respect to $t$ is given by : \bq\label{log1}\dfrac{H^{'}(x,-\frac{t^2}{4})}{H(x,-\frac{t^2}{4})}=-\sum_{j=1}^n\dfrac{\frac{1}{2}x_jt}{1+\frac{1}{4}x_jt^2}=i\sum_{m=1}^\infty p_{m}(x)\left(\dfrac{it}{2}\right)^{2m-1}.\eq 

\bigskip

Let $\Gamma$ be the algebra of symmetric functions on $${\Bbb C}^{(\infty)}=\{z_1,z_2,\dots)\, |\, z_i\in\Bbb C \ \textrm{are zero for {\it i} large enough}\}.$$ The set $\{s_{\bf m}\}$ where ${\bf m}$ runs over all partitions is a system of linear generators of $\Gamma$. Furthermore, the sets $\{h_m\}_{m\geq1}$ and $\{p_m\}_{m\geq1}$ are systems of generators of $\Gamma$ : $$\Gamma=\Bbb C [h_1,h_2,\dots]=\Bbb C [p_1,p_2,\dots].$$

 Let us consider the algebra morphism : \bq\label{log2}\Gamma \rightarrow \mathscr{C}(\Omega),\, g\mapsto\widetilde{g},\eq which is uniquely determined by : \bq\label{log3}\widetilde{p}_1(\omega)=\gamma+\sum_{j=1}^\infty\alpha_j, \quad \quad \widetilde{p}_m(\omega)=\sum_{j=1}^\infty\alpha_j^m \quad (m\geq2).\eq

The functions $\widetilde{p}_m$ are continuous on $\Omega$. This can be shown by taking in (\ref{log}) $f\equiv1$ for $m=1$, or $f(t)=t^{2m-2}$ for $m\geq2$.

\bigskip

\begin{pro}\label{nes1} It holds that $$\widetilde{H}(\omega,-\frac{\lambda^2}{4})=\Pi(\omega,\lambda).$$ Furthermore, the Taylor expansion of $\Pi(\omega,\lambda) $ is : $$\Pi(\omega,\lambda)=\sum_{m=0}^\infty \widetilde{h}_m(\omega)\left(-\frac{\lambda^2}{4}\right)^m.$$ Here $\widetilde{H}(\omega,-\frac{\lambda^2}{4})$ is the image of $H(.,-\frac{\lambda^2}{4})$ under the morphism {\rm(\ref{log2})} and $\Pi(\omega,\lambda)$ is the modified P{\'o}lya function.\end{pro}
\begin{paragraph}{\bfseries{Proof.}} {\rm  By (\ref{logg}), (\ref{log3}) and (\ref{log1}), we have $$\dfrac{\Pi^{'}(\omega,\lambda)}{\Pi(\omega,\lambda)}=i\sum_{m=1}^\infty\widetilde{p}_{m}(\omega)\left(\dfrac{i\lambda}{2}\right)^{2m-1}=\dfrac{\widetilde{H}^{'}(\omega,-\frac{\lambda^2}{4})}{\widetilde{H}(\omega,-\frac{\lambda^2}{4})}.$$ Since $\Pi(\omega,0)=1$ and $\widetilde{H}(\omega,0)=1$, the statement follows.$\hspace{1cm}\Box$\\}
\end{paragraph}

Let us consider now the map $$T_n\; : {\Bbb R}^n\rightarrow\Omega;\; \big(\lambda_1,\lambda_2,\dots,\lambda_n\big)\mapsto\omega=(\alpha,\gamma)$$ given by

$$\alpha_j=\left(\frac{\lambda_j}{n}\right)^2,\quad \quad \gamma=0.$$

\begin{teo}\label{nes2}Let $\lambda^{(n)}\in\Bbb R^n$ be a sequence such that the following limit exists for the topology of $\Omega$ : $$\lim_{n\to\infty}T_n\left(\lambda^{(n)}\right)=\omega .$$ Then, for every $g\in\Gamma$, homogeneous of degree $m$, $$\lim_{n\to\infty}\dfrac{1}{n^{2m}}\;g\left(\left(\lambda^{(n)}\right)^2\right)=\widetilde{g}(\omega).$$ 
\end{teo}
\begin{paragraph}{\bfseries{Proof.}} {\rm It is enough to prove the result for $g=p_m$ since the Newton power sums generate $\Gamma$. Let $m=1$, then $$p_1(\lambda^{(n)})=\sum_{j=1}^n\lambda_j^{(n)}, \quad \quad \widetilde{p}_1(\omega)=\gamma+\sum_{j=1}^\infty\alpha_j.$$  By assumption, for every continuous function $f$ on $\Bbb R$, $$\lim_{n\to\infty}\sum_{k=1}^n\left(\dfrac{\lambda_k^{(n)}}{n}\right)^2 f\left(\dfrac{\lambda_k^{(n)}}{n}\right)= \gamma f(0)+\sum_{k=1}^\infty \alpha_k f(\alpha_k) .$$ In particular, by taking $f\equiv1$ we get  $$\lim_{n\to\infty}\sum_{k=1}^n\left(\dfrac{\lambda_k^{(n)}}{n}\right)^2= \gamma+\sum_{k=1}^\infty \alpha_k,$$ $$\lim_{n\to\infty}\dfrac{1}{n^2}\;p_1\left(\left(\lambda^{(n)}\right)^2\right)=\widetilde{p}_1(\omega) .$$    For $m\geq2$, by taking $f(t)=t^{2m-2}$, one obtains $$\lim_{n\to\infty}\dfrac{1}{n^{2m}}\;p_m\left(\left(\lambda^{(n)}\right)^2\right)=\widetilde{p_m}(\omega).\hspace{1cm}\Box$$\\}

\end{paragraph}

\begin{pro}\label{nes4} For any $\omega\in\Omega$ and $\xi_1,\dots,\xi_k\in\Bbb C$, it holds that $$\sum_{{\bf m}:{\rm partitions}}\widetilde{s}_{\bf m}(\omega)s_{\bf m}\left(-\dfrac{\xi_1^2}{4},\dots,-\dfrac{\xi_k^2}{4}\right)=\prod_{j=1}^k\Pi(\omega,\xi_j).$$  \end{pro}
\begin{paragraph}{\bfseries{Proof.}} {\rm  Recall the Cauchy formula (see \cite{Far2}, Proposition 2.5) : $$\sum_{\bf m}s_{\bf m}(x_1,x_2,\dots)s_{\bf m}(y_1,\dots,y_k)=\prod_{i=1}^\infty\prod_{j=1}^k\frac{1}{1-x_iy_j}=\prod_{j=1}^k H(x,y_j).$$ We apply the morphism (\ref{log2}) to both sides of the preceding equality with $y_j=-\dfrac{\xi_j^2}{4}$. The result follows by Proposition \ref{nes1}.$\hspace{1cm}\Box$\\        }
\end{paragraph}

\bigskip

\begin{teo}\label{nes5} As in Theorem \ref{nes2}, we assume $$\lim_{n\to\infty}T_n\left(\lambda^{(n)}\right)=\omega .$$ Then, for a fixed diagonal matrix $\xi={\rm diag}(\xi_1,\dots,\xi_k,0,\dots)$, the sequence $\pp_n(\xi,\lambda^{(n)})$ converges uniformly on compact sets in $\Bbb R^k$ : $$\lim_{n\to\infty}\pp_n(\xi,\lambda^{(n)})=\prod_{j=1}^k\Pi(\omega,\xi_j).$$

\end{teo}
\begin{paragraph}{\bfseries{Proof.}} {\rm  Consider first $k=1$. This case corresponds to a single variable $\xi=(\xi,0,\dots)$ and the unique non-zero terms in the Schur function expansion are those for which ${\bf m}=(m,0,\ldots,0)$. Hence by (\ref{nes3}) we have $$\pp_n(\xi,\lambda^{(n)})=\sum_{m=0}^\infty \left(\frac{(n-1)!}{(m+n-1)!}\right)^2\ h_m({\lambda^{(n)}}^2)\left(-\dfrac{\xi^2}{4}\right)^m.$$  Since $$\dfrac{(n-1)!}{(m+n-1)!}\sim\dfrac{1}{n^m} \quad {\rm as}\ n\to\infty ,$$ by Theorem \ref{nes2} $$\lim_{n\to\infty}\left(\frac{(n-1)!}{(m+n-1)!}\right)^2\ h_m({\lambda^{(n)}}^2)=\widetilde{h}_m(\omega).$$ Now by applying Proposition \ref{proofc} about the convergence of $\mathcal C^\infty$-functions of positive type on $\Bbb R^d$, we obtain  $$\lim_{n\to\infty}\pp_n(\xi,\lambda^{(n)})=\sum_{m=0}^\infty \widetilde{h}_m(\omega)\left(-\dfrac{\xi^2}{4}\right)^m.$$ Finally, by Proposition \ref{nes1}, $$\sum_{m=0}^\infty \widetilde{h}_m(\omega)\left(-\frac{\lambda^2}{4}\right)^m=\Pi(\omega,\lambda).$$ Now, let us consider the multivariable case $\xi=(\xi_1,\dots,\xi_k,0,\dots)$ for which $k>1$. If $m_{k+1}> 0$, then $s_{\bf m}(\xi_1,\dots,\xi_k,0,\dots)=0$. In consequence $$\pp_n(\xi,\lambda^{(n)})=\sum_{m_1\geq\dots\geq m_k\geq0}\left(\dfrac{{\bf \delta!}}{(\bf m+\bf \delta)!}\right)^2s_{\bf m}\left({\lambda^{(n)}}^2\right)s_{\bf m}\left(-\dfrac{\xi_1^2}{4},\dots,-\dfrac{\xi_k^2}{4},0,\dots\right).$$ But, for ${\bf m}$ fixed, $n\to\infty$,  $$\dfrac{{\bf \delta!}}{(\bf m+\bf \delta)!}\sim\dfrac{1}{n^{|{\bf m}|}},\quad |{\bf m}|=m_1+m_2+\dots.$$ Hence by Theorem \ref{nes2} $$\lim_{n\to\infty}\left(\dfrac{{\bf \delta!}}{(\bf m+\bf \delta)!}\right)^2s_{\bf m}\left({\lambda^{(n)}}^2\right)=\widetilde{s}_{\bf m}(\omega).$$ Similarly, by Proposition \ref{proofc},   $$\lim_{n\to\infty}\pp_n(\xi,\lambda^{(n)})=\sum_{{\bf m}:{\rm partitions}}\widetilde{s}_{\bf m}(\omega)s_{\bf m}\left(-\dfrac{\xi_1^2}{4},\dots,-\dfrac{\xi_k^2}{4},0,\dots\right).$$ Finally, by Proposition \ref{nes4}, $$\sum_{{\bf m}:{\rm partitions}}\widetilde{s}_{\bf m}(\omega)s_{\bf m}\left(-\dfrac{\xi_1^2}{4},\dots,-\dfrac{\xi_k^2}{4},0,\dots\right)=\prod_{j=1}^k\Pi(\omega,\xi_j).\hspace{1cm}\Box$$      }
\end{paragraph}

\bigskip

The preceding theorem shows that the limit of a spherical function of positive type on $V_n$ is a spherical function of positive type on $V_\infty$ which is given as a finite product of modified P{\'o}lya functions. In order to prove that all spherical functions of positive type on $V_\infty$ are obtained in the same way, we need to prove the converse of Theorem \ref{nes5}. This will be done in Theorem \ref{bi} using the following lemma :\\

\begin{lem}\label{keyrev} For every $n\geq1$, let $\omega^{(n)}$ be the point in $\Omega$ associated to an orbital measure $\n^{(n)}$ that weakly converges to an ergodic measure $\nu$ on $V^\infty$. If \bq\label{big}\sup_{n} \Big(p_1(\alpha^{(n)})\Big)=R<\infty,\eq then the sequence of parameters $\omega^{(n)}$ converges in $\Omega$.\\
\end{lem}

\begin{paragraph}{\bfseries{Proof.}}{\rm By condition (\ref{big}), the sequence $\omega^{(n)}$ belongs to the compact set $\Omega_R$ (see Corollary \ref{compact}). We can then extract a subsequence $\left(\omega^{(n_k)}\right)_k$ which converges in $\Omega$ to $\omega=(\alpha,\gamma)$. Hence, by Theorem \ref{nes5}, the Fourier transform $\pp^{(n_k)}$ of $\nu^{(n_k)}$ uniformally converges on compact sets : $$\lim_{k\rightarrow\infty}\pp^{(n_k)}\left({\rm diag}\left(\xi,0,\dots\right)\right)=\Pi(\omega,\xi).$$ Therefore, the Fourier transform $\pp$ of $\nu$ is identically equal to $\Pi(\omega,.).$\\

If another sub-sequence such that $\omega^{({n_k}^{'})}$ converges in $\Omega$ to $\omega^{'}$, then $\pp\equiv\Pi(\omega^{'},.)$. By uniqueness, we get $\omega=\omega^{'}$, and $\omega^{(n)}$ do have in consequence a unique accumulation point. Therefore, it necessarily converges to $\omega$.$\hspace{1cm}\Box$\\}
\end{paragraph}

\begin{teo}\label{bi}
The spherical functions of positive type on $V_\infty$ (i.e. the extreme points of $\mathfrak P$) are the functions $\varphi_\omega$ defined, for every $k\geq1$, by : $$\pp_\omega\big({\rm diag}(\xi_1,\dots,\xi_k,0,\dots)\big)=\Pi(\omega,\xi_1)\dots\Pi(\omega,\xi_k),$$  where $\Pi(\omega,.)$ is the modified P{\'o}lya function associated to $\omega\in\Omega$.\\
\end{teo}

\begin{paragraph}{\bfseries{Proof.}}{\rm
(a) The function $\pp_{\omega}$, which is of positive type, is spherical since it is multiplicative (Theorem \ref{multi}). Therefore, it is an extreme point in $\mathfrak P$.

(b) Let $\varphi\in {\rm ext}(\mathfrak P)$. It is the Fourier transform of an ergodic measure $\nu$ on $V^\infty$ relatively to the action of $K_\infty$. By using a theorem due to A. Vershik (\cite{Ol2}, Theorem 3.2), the measure $\nu$ is the weak limit of $\nu^{(n)}$, where $\nu^{(n)}$ is a sequence of orbital measures relatively to $K_n$. Hence $\varphi$
 is the uniform limit on compact sets of the sequence $\pp_n$, where $\pp_n$ is the Fourier transform of $\nu^{(n)}.$ In particular we have $$\lim_{n\rightarrow\infty}\pp_n\left({\rm diag}\left(\xi,0,\dots\right),\lambda^{(n)}\right)=\pp\left({\rm diag}\left(\xi,0,\dots\right)\right).$$ We will prove that $\omega^{(n)}=T_n(\lambda^{(n)})$ converges in $\Omega$, which gives the converse of Theorem \ref{nes5}. Let us suppose now that the condition (\ref{big}) is not satisfied :

 $$\sup_n\left(p_1\left(\left(\dfrac{\lambda^{(n)}}{n}\right)^2\right)\right)=\infty.$$ 

There exists a positive sequence $\ep_n$ such that $$\lim_{n\rightarrow\infty}\ep_n=0\ {\rm and}\ \lim_{n\rightarrow\infty}p_1\left(\ep_n\left(\dfrac{\lambda^{(n)}}{n}\right)^2\right)=1.$$ One can remark that multiplying $(\dfrac{\lambda^{(n)}}{n})^2$ by $\ep_n$ is the same as multiplying $\xi$ by $\ep_n$. Hence, by Lemma \ref{keyrev}, there exists $\omega\in\Omega$ such that :

 $$\lim_{n\rightarrow\infty}\pp_n\left({\rm diag}\left(\ep_n\xi,0,\dots\right),\lambda^{(n)}\right)=\Pi(\omega,\xi).$$

 The modified P{\'o}lya function $\Pi(\omega,.)$ is not identically equal to $1$, because $\omega\ne0$. Therefore, in a neighborhood $D(0,R)$ of $0$, the function $\Pi(\omega,.)$ is not identically equal to $1$. But this leads to a contradiction, because, by our assumption, \ba\lefteqn{\lim_{n\rightarrow\infty}\pp_n\left({\rm diag}\left(\ep_n\xi,0,\dots\right),\lambda^{(n)}\right)}\\& &=\lim_{n\rightarrow\infty}\sum_{m=0}^\infty \left(\frac{n^m\ep_n^m}{n(n+1)\dots(n+m-1)}\right)^2h_m\left(\left(\dfrac{\lambda^{(n)}}{n}\right)^2\right)\;\left(-\dfrac{\xi^2}{4}\right)^m=1,\ea and then $\Pi(\omega,.)\equiv1.$ Hence, there exists $\omega\in\Omega$ such that : $\pp(\xi)=\pp_{\omega}(\xi)\ \ (\xi\in V_\infty).\hspace{1cm}\Box$}
\end{paragraph}

\bigskip

\begin{nota} {\rm  A measurable space is said to be standard if it is isomorphic to a Borel subset in a polish space which is equipped with the $\sigma$-algebra induced by the Borel one. One can prove that the correspondence $\Omega\leftrightarrow$ {\rm ext(}$\mathfrak P${\rm)} is an isomorphism between two standard spaces. This enables us to prove a parameterized version of the generalized Bochner theorem (\cite{Rab}, Theorem 7) : let $\varphi$ be a $K_\infty$-invariant continuous function of positive type on $V_\infty$ with $\varphi(e)=1$. Then, there exists a unique probability measure $\m$ defined on $\Omega$ such that, for every $g\in V_\infty$,
$$\varphi(g)=\int_{\Omega}\varphi_{\omega}(g)\, \m(d\omega).$$}

\end{nota}

\begin{paragraph}{\bfseries {Acknowledgement}}{\rm I would like to express my gratitude to Professor Jacques Faraut for numerous suggestions during the preparation of this paper.}
\end{paragraph}

\bibliographystyle{amsalpha}

\end{document}